\documentclass[a4paper,12pt]{article}
\usepackage{fancyhdr,graphicx}
\usepackage{amsfonts}
\usepackage{amssymb}
\usepackage{amsthm}
\usepackage{newlfont}
\usepackage{amsmath}
\usepackage[top=2.5cm,bottom=3.0cm,left=2.5cm,right=2.5cm]{geometry}
\usepackage{multirow}
\usepackage{array}
\usepackage{longtable}
\usepackage{bm}
\usepackage{setspace}

\newcommand{\trou}{\vspace{1 mm}}
\newcommand{\noi}{\noindent}

\renewcommand{\baselinestretch}{1.4}

\usepackage{latexsym,bm}
\title{\vspace{0cm}{The index of a string consisting of 4 blocks}}

\author{Jianxin Wei $^{a}$, Heping Zhang $^{b}$  \\
\footnotesize{$^{a}$ School of Mathematics and Statistics Science, Ludong University}\\
\footnotesize{Yantai, Shandong 264025, P. R. China}\\
\footnotesize {$^{b}$ School of Mathematics and Statistics, Lanzhou University}\\
 \footnotesize {Lanzhou, Gansu 730000, P. R. China}\\
\small{E-mail addresses: wjx0426@163.com, zhanghp@lzu.edu.cn}}
\date{}

\begin{document}

\maketitle

\begin{abstract}
Generalized Fibonacci cube $Q_{d}(f)$, introduced by Ili\'{c}, Klav\v{z}ar and Rho,
is the graph obtained from the $d$-hypercube $Q_{d}$ by removing all vertices that contain $f$ as a substring.
The smallest integer $d$ such that $Q_{d}(f)$ is not
an isometric subgraph of $Q_{d}$ is called the index of $f$.
A non-extendable sequence of contiguous equal digits in a string $\mu$ is called a block of $\mu$.
The question that determine the index of a string consisting of at most 3 blocks is solved
by Ili\'{c}, Klav\v{z}ar and Rho.
This question is further studied and the index of a string consisting of 4 blocks is determined,
and the necessity of a string being good is also given for the strings with even blocks.
\end{abstract}
\textbf{Key words:} Generalized Fibonacci cube; Index of a string; Good string; Bad string.

\section{Introduction and Preliminaries}

Hsu \cite{whs} introduced \emph{Fibonacci cube} as a model for interconnection networks,
which has similar properties as hypercube.
The vertex set of Fibonacci cube $\Gamma_{d}$ is the set of all binary strings $b_{1}b_{2}\ldots b_{d}$ containing no two consecutive $1$s
and two vertices are adjacent in $\Gamma_{d}$ if they differ in precisely one bit.
Klav\v{z}ar and \v{Z}igert \cite{sp} applied Fibonacci cubes in chemical graph theory and showed those cubes are precisely the resonance graphs of fibonaccenes.
More generally Zhang et al. \cite{zoy} described the class of planar bipartite graphs that have Fibonacci cubes as their resonance graphs.
For more about Fibonacci cubes, see \cite{sk} for a survey.

A binary string $f$ is called a \emph{factor} of binary string $\mu$ if $f$ appears as a
sequence of $|f|$ consecutive bits of $\mu$, where $|f|$ denotes the length of $f$.
$\Gamma_{d}$ can be seen as the graph obtained from
$Q_{d}$ by removing all strings that contain $11$ as a factor.
Inspired by this, Ili\'{c}, Klav\v{z}ar and Rho \cite{asy} introduced \emph{generalized Fibonacci cube}, $Q_{d}(f)$, as the graph obtained from $Q_{d}$ by removing all strings that contain $f$ as a factor, where $f$ is some given binary string.
In this notation Fibonacci cube $\Gamma_{d}$ is the graph $Q_{d}(11)$.
The subclass $Q_{d}(1^{s})$ of generalized Fibonacci cube has been studied in \cite{liu,nzs}.

For a connected graph $G$, the \emph{distance} $d_{G}(\mu, \nu)$ between vertices $\mu$ and $\nu$ is the length of a shortest $\mu,\nu$-path.
Given two binary strings $\alpha$ and $\beta$ with the same length,
their \emph{Hamming distance} $H(\alpha,\beta)$ is the number of bits in which they differ.
It is known that \cite{ik} for any vertices $\alpha$ and $\beta$ of $Q_{d}$, $H(\alpha,\beta)=d_{Q_{d}}(\alpha,\beta)$.

Obviously, for any subgraph $H$ of $G$, $d_{H}(\mu,\nu)\geq d_{G}(\mu,\nu)$.
If $d_{H}(\mu,\nu)=d_{G}(\mu,\nu)$ for all $\mu,\nu\in V(H)$,
then $H$ is called \emph{an isometric subgraph} of $G$,
and simply write $H\hookrightarrow G$, and $H\not\hookrightarrow G$ otherwise.

A string $f$ is \emph{good} if $Q_{d}(f)\hookrightarrow Q_{d}$ for all $d\geq1$,
it is \emph{bad} otherwise \cite{sks}.
The \emph{index} of a binary string $f$, denoted $B(f)$,
is the smallest integer $d$ such that $Q_{d}(f)\not\hookrightarrow Q_{d}$ \cite{asy}.
Obviously, $f$ is good if and only if $B(f)=+\infty$ and $f$ is bad if and only if $B(f)<+\infty$.
It was shown that for about eight percent of all strings are good \cite{sks}.

A non-extendable sequence of contiguous equal digits
in a string $\alpha$ is called a \emph{block} of $\alpha$.
Let $F'=\{\ f \ | \ f$ is a string consisting of at most 3 blocks$\}$.
The question that determine the index of $f\in F'$ was solved by Ili\'{c},
Klav\v{z}ar and Rho \cite{asy}. The result is shown in Table 1.

{\renewcommand\baselinestretch{0.9}\selectfont
\begin{center} \fontsize{10.5pt}{11.5pt}\selectfont
\begin{longtable}{l|l|l|l|l}
\caption{\label{comparison} Classification of the index of the string $f\in F'$.}\\
\hline
$(i')$ & $r$ & $s$ & $t$ & $B(f)$\\
\hline
$(1')$ & $r\geq1$ & $s=0$ & $t=0$ &\multirow{2}{*} {$+\infty$} \\
\cline{1-4}
$(2')$ & $r\geq1$ & $s=1$ & $t=0$ & \\
\hline
$(3')$& $r=2$ & $s\geq2$ & $t=0$ & $s+5$ \\
\hline
$(4')$& $r\geq3$ & $s\geq3$ & $t=0$ & $2r+2s-2$ \\
\hline
$(5')$& $r\geq1$ & $s\geq1$ & $t\geq1$ & $r+s+t+1$ \\
\hline
\end{longtable}
\end{center}
\par}

This question is further studied for the strings consisting of even blocks in this paper.
Let $f=1^{x_{1}}0^{y_{1}}1^{x_{2}}0^{y_{2}}\cdots1^{x_{n}}0^{y_{n}}$,
where $x_{i}\geq 1$, $y_{i}\geq 1$, $i=1,\ldots,n$ and $n\geq2$.
We find 5 classes of bad strings and give a necessity of a string being good in the following theorem.

\trou\noi{\bf Theorem 1.1.}
\emph{If $f$ is good, then it satisfies one of the following cases:}

(a) \emph{$x_{1}=1$ and $y_{n}=y_{1}+1$},

(b) \emph{$x_{1}=x_{n}=1$ and $y_{1}=y_{n}$},

(c) \emph{$x_{1}=1$, $x_{n}\geq2$ and $y_{1}=y_{n}+1$} and

(d) \emph{$x_{1}=1$, $x_{n}=2$ and $y_{1}\geq y_{n}+2$.}

We pay special attention to the strings consisting of 4 blocks.
Let $F=\{f=1^{r}0^{s}1^{t}0^{k}|$ $r$, $s$, $t$ and $k\geq1 \}$.
We get the following result.

\trou\noi{\bf Theorem 1.2.}
\emph{Let $f\in F$. Then $B(f)$ is given just as shown in Table 2.}

{\renewcommand\baselinestretch{0.9}\selectfont
\begin{center} \fontsize{10.5pt}{11.5pt}\selectfont
\begin{longtable}{l|l|l|l|l|l}
\caption{\label{comparison} Classification of the index of the string $f\in F$.}\\ 
\hline
(i) & $r$ & $t$ & $s$ &$k$ & $B(f)$\\
\hline
(1)  & $r=1$ & $t\geq 1$ & $s\geq 1$ & $k=s+1$ & \multirow{4}{*} {$+\infty$} \\
\cline{1-5}
(2) & $r=1$ & $t=1$ & $s=k$ & $k\geq1$ &  \\
\cline{1-5}
(3) & $r=1$ & $t\geq 2$ & $s=k+1$ & $k\geq1$ &  \\
\cline{1-5}
(4) & $r=1$ & $t=2$ & $s\geq k+2$ & $k\geq2$ & \\
\hline
(5) & $r\geq t+3$ & $t\geq 1$ & $s\geq3$ & $k=2$ & $ 2r+2s+t+2$ \\
\hline
(6) & $r=t+2$ & $t\geq 1$ &  $s\geq 2$  & $s\geq k\geq 1$ & $  3t+2s+k+2$ \\
\hline
\multirow{3}{*} {(7)} & $r\geq t+3$ & $t\geq 1$ & $s\geq 3$ & $k=1$ &
\multirow{3}{*} {$ 2r+2s+t+k-1$}\\
\cline{2-5} & $r=t+2$ & $t\geq 1$  & $s=1$ & $k=1$ &  \\
\cline{2-5} &$r=t+2$ & $t\geq 1$  & $s\geq 2$ & $k=s+1$ &  \\
\hline
\multirow{3}{*} {(8)} & $r= t+1$ & $t\geq 2$ & $s\geq 1$ & $k=s+1$ &
 \multirow{3}{*} {$ r+2s+2t+k$} \\
\cline{2-5} & $r=2$ & $t=1$ & $s\geq 1$ & $k\geq s+1$ &\\
\cline{2-5} & $r\geq 3$ & $t=1$ & $s=1$ & $k\geq 3$ &\\
\hline
(9) & $r=1$ & $t=1$ & $s\geq 1$ & $k\geq s+3$ & $ 2s+k+4$ \\
\hline
(10) & $r=1$ & $t=1$ & $s\geq k+1$ & $k\geq 1$ & $2k+s+3$ \\
\hline
(11) & $r\geq 2$ & $t=r$ & $s\geq k+1$ & $k\geq 1$ &$ 3r+2s+k-1$ \\
\hline
(12) & $r\geq 1$ & $t\geq r+2$ & $s\geq k+2$ & $k\geq 1$ & $ 2s+2t+r+k-2$ \\
\hline
(13)& $2\geq r\geq 1$ & $t=2$ & $s\geq 1$ & $k\geq s+3$ & $ 2s+k+r+4$ \\
\hline
\multirow{4}{*}{(14)}  & $r\geq 3$  & $t=2$ & $s\geq 1$ &
\multirow{4}{*} {$k\geq s+3$} & \multirow{5}{*} {$ 2(r+s+t+k-1)$}\\
\cline{2-4} & $r=3$  & $t=1$ & $s\geq2$ & \\
\cline{2-4} & $r\geq t$ & $t\geq 3$ & $s\geq 1$  &\\
\cline{2-4} & $r\geq 3$ & $t\geq r+1$ & $s\geq 1$ & \\
\hline
\multirow{4}{*}{(15)} & $r\geq 2$ & $t=r+1$ & $s\geq k+1$ & $k\geq 2$ &
\multirow{4}{*} {$ r+2s+2t+k+2$}\\
\cline{2-5} & $r=t+1$ & $t\geq 1$ & $s\geq k$ & $k\geq 2$  \\
\cline{2-5} & $r=t+2$ & $t\geq 1$ & $s\geq 2$ & $k=s+2$ \\
\cline{2-5} & $r=t $ & $t\geq 2$ & $s\geq 2$ & $k=s$  \\
\hline
\end{longtable}
\end{center}
\par}

In the rest of this section some necessary definitions and results are introduced.
With $e_{i}$ we denote the binary string with $1$ in the $i$-th bit and 0 elsewhere.
For strings $\alpha$ and $\beta$ of the same length
let $\alpha+\beta$ denote their sum computed bitwise modulo 2.
In particular, $\alpha+e_{i}$ is the string obtained from $\alpha$ by reversing its $i$-th bit.
The null string, denoted by $\lambda$, is a string of length zero.
For the convenience of use,
let $\mu_{s,t}$ be the factor of $\mu$ that starts from the $s$-th bit to the $t$-th bit if $t\geq s$
and $\lambda$ if $t<s$, where $s\leq|\mu|$ and $t\leq |\mu|$.

Let $\alpha =a_{1}a_{2}\cdots a_{d_{1}}$ and $\beta=b_{1}b_{2}\cdots b_{d_{2}}$,
with $\alpha\beta$ we denote $a_{1}a_{2}\cdots a_{d_{1}}b_{1}b_{2}\cdots b_{d_{2}}$,
with $\alpha^{R}=a_{d_{1}}\cdots a_{2}a_{1}$ denote the reverse of $\alpha$ and $\overline{\alpha}=\overline{a_{1}}\ \overline{a_{2}} \cdots \overline{a_{d_{1}}}$ the complement of $f$, where $\overline{f_{i}}=1-f_{i}$, $i=1,\ldots,d$.

\trou\noi{\bf Proposition 1.3 (\cite{asy}).} \emph{Let f be any string and $d\geq 1$.
Then $Q_{d}(f)\cong Q_{d}(\overline{f})\cong Q_{d}(f^{R})$.}

For two vertices $\mu$ and $\nu$ of graph $G$,
the set of vertices lying on shortest $\mu,\nu$-paths is called the interval between $\mu $
and $\nu$, denoted by $I_{G}(\mu,\nu)$.
Let $\alpha$ and $\beta\in V(Q_{d}(f))$ and $p\geq2$.
Then $\alpha$ and $\beta$ are called $p$-critical words \cite{asy} for $Q_{d}(f)$ if $d_{Q_{d}}(\alpha,\beta)=H(\alpha,\beta)=p$,
but none of the neighbors of $\alpha$ in $I_{Q_{d}}(\alpha,\beta)$
belongs to $Q_{d}(f)$ or none of the neighbors of $\beta$ in
$I_{Q_{d}}(\alpha,\beta)$ belongs to $Q_{d}(f)$.
The following proposition gives a tool to prove $Q_{d}(f)\not\hookrightarrow Q_{d}$.

\trou\noi{\bf Proposition 1.4 (\cite{asy}).}
\emph{If there exist $p$-critical words for $Q_{d}(f)$ for some $p\geq2$, then $Q_{d}(f)\not\hookrightarrow Q_{d}$.}

\trou\noi{\bf Proposition 1.5 (\cite{sks}).} \emph{Suppose that $Q_{d}(f )\not\hookrightarrow Q_{d}$. Then $Q_{d'}(f)\not\hookrightarrow Q_{d'}$ for all $d'>d$.}

Note that if $f$ is bad, then $Q_{d}(f)\hookrightarrow Q_{d}$ only for $d<B(f)$ by Proposition 1.5.

We proceed as follows.
In the next section, we study the strings with even blocks.
A few special classes of bad strings are given and
Theorem 1.1 is proved.
In the last section, Theorem 1.2 is proved,
in other words, the index of  every string consisting of 4 blocks is determined.

\section{Strings consisting of $2n$ blocks}

In this section we give 5 classes of bad strings with $2n$ blocks,
which will be used to give the necessity of a string being good,
and also to classify the index of a string from $F$ in the next section.
All the strings considered in the following lemmas and theorem consist of $2n$ blocks.
Assume without loss of generality that $f=1^{x_{1}}0^{y_{1}}1^{x_{2}}0^{y_{2}}\cdots1^{x_{n}}0^{y_{n}}$ by Proposition 1.3,
where $x_{i}\geq 1$, $y_{i}\geq 1$, $i=1,\ldots,n$ and $n\geq2$.

\trou\noi{\bf Lemma 2.1.}
\emph{If $x_{1}\geq2$ and $y_{n}\geq2$,
then $Q_{d}(f)\not\hookrightarrow Q_{d}$ for
$d\geq 2\sum\limits_{i=1}^{n}(x_{i}+y_{i})-2$.}

\trou\noi{\bf Proof.} Let $d_{0}=2\sum\limits_{i=1}^{n}(x_{i}+y_{i})-2$,

$\alpha=1^{x_{1}}0^{y_{1}}1^{x_{2}}0^{y_{2}}\cdots1^{x_{n}} 0^{y_{n}-2}101^{x_{1}-2}0^{y_{1}}1^{x_{2}}0^{y_{2}}\cdots1^{x_{n}}0^{y_{n}}$ and

$\beta=1^{x_{1}}0^{y_{1}}1^{x_{2}}0^{y_{2}}\cdots1^{x_{n}}0^{y_{n}-2}01
1^{x_{1}-2}0^{y_{1}}1^{x_{2}}0^{y_{2}}\cdots1^{x_{n}}0^{y_{n}}$.

Note that $|\alpha|=|\beta|=d_{0}$, $H(\alpha,\beta)=2$
and the only vertices on the two shortest $\alpha,\beta$-paths in $Q_{d_{0}}$ are

$\mu=1^{x_{1}}0^{y_{1}}1^{x_{2}}0^{y_{2}}\cdots1^{x_{n}}0^{y_{n}}
1^{x_{1}-2}0^{y_{1}} 1^{x_{2}}0^{y_{2}}\cdots1^{x_{n}}0^{y_{n}}$ and

$\nu=1^{x_{1}}0^{y_{1}}1^{x_{2}}0^{y_{2}}\cdots1^{x_{n}}0^{y_{n}-2}1^{x_{1}}$
$0^{y_{1}}1^{x_{2}}0^{y_{2}}\cdots1^{x_{n}}0^{y_{n}}$.

But none of $\mu$ and $\nu$ is a vertex of $Q_{d_{0}}(f)$.
Now we claim that $\alpha$ and $\beta$ are vertices of $Q_{d_{0}}(f)$.
If this claim holds, then $\alpha$ and $\beta$ are 2-critical words for $Q_{d_{0}}(f)$,
hence $Q_{d_{0}}(f)\not\hookrightarrow Q_{d_{0}}$ by Proposition $1.4$.
To prove this claim we can check every factor of consecutive $2n$ blocks
contained in $\alpha$ or $\beta$ is $f$ (or not) by a system of equations and inequalities has a solution of positive integers (or not).
Obviously, we only need consider the factors beginning the block consisting of 1s.
There are four cases to be considered:
(a) $x_{1}\geq3$ and $y_{n}\geq3$;
(b) $x_{1}\geq3$ and $y_{n}=2$;
(c) $x_{1}=2$ and $y_{n}\geq3$ and
(d) $x_{1}=2$ and $y_{n}=2$.

Since cases (b), (c) and (d) can be proved by similar arguments as in case (a),
we only give the proof of case (a) to save space.

There are $n+2$ systems of equations and inequalities for $\alpha$. The first and second are:

{\renewcommand\baselinestretch{0.8}\selectfont
\begin{equation}\nonumber
(1) \begin{cases}
x_{t}\leq x_{1}, t=1,\\
x_{t}=x_{t},t=2,\ldots,n,\\
y_{t}=y_{t},t=1,\ldots,n-1,\\
y_{t}\leq y_{n}-2,t=n.\\
\end{cases}
(2) \begin{cases}
x_{t}\leq x_{2}, t=1,\\
x_{t}=x_{t+1},t=2,\ldots,n-1,\\
x_{t}=1,t=n,\\
y_{t}=y_{t+1},t=1,\ldots,n-2,\\
y_{t}=y_{n}-2,t=n-1,\\
y_{t}\leq 1,t=n.\\
\end{cases}
\end{equation}

The $i$-th $(i=3 ,\ldots, n,n\geq3)$ system of equations and
 inequalities is:

{\renewcommand\baselinestretch{0.8}\selectfont
\begin{equation}\nonumber
(i) \begin{cases}
x_{t}\leq x_{i}, t=1,\\
x_{t}=x_{t+i-1},t=2,\ldots,n-i+1,\\
x_{t}=1,t=n-i+2,\\
x_{t}=x_{1}-2,t=n-i+3,\\
x_{t}=x_{t+i-n-2},t=n-i+4,\ldots,n,\\
y_{t}=y_{t+i-1},t=1,\ldots,n-i,\\
y_{t}=y_{n}-2,t=n-i+1,\\
y_{t}=1,t=n-i+2,\\
y_{t}=y_{t+i-n-2},t=n-i+3,\ldots,n-1,\\
y_{t}\leq y_{i-2},t=n.\\
\end{cases}
\end{equation}
\par}

The last two are:

\begin{equation}\nonumber
(n+1) \begin{cases}
x_{t}\leq 1, t=1,\\
x_{t}=x_{1}-2,t=2,\\
x_{t}=x_{t-1},t=3,\ldots,n,\\
y_{t}=1,t=1,\\
y_{t}=y_{t-1},t=2,\ldots,n-1,\\
y_{t}\leq y_{t-1},t=n.\\
\end{cases}
(n+2) \begin{cases}
x_{t}\leq x_{1}-2, t=1,\\
x_{t}=x_{t},t=2,\ldots,n,\\
 y_{t}=y_{t-1},t=1,\ldots,n-1,\\
y_{t}\leq y_{t},t=n.\\
\end{cases}
\end{equation}

\par}

In the first system,
the last inequality is $y_{n}\leq y_{n}-2$,
and in the $(n+2)$-th system, the first inequality is $x_{1}\leq x_{1}-2$,
obviously those are impossible.
If we add up the first to $n$-th inequalities and equations in the $i$-th system,
then $x_{i}\leq-1$, $i=3,\ldots,n$. It is impossible.
Hence all those systems have no solution of positive integers,
that is, $\alpha\in V(Q_{d_{0}}(f))$.
It can be shown that $\beta\in V(Q_{d_{0}}(f))$ similarly.
So $\alpha$ and $\beta$ are 2-critical words for $Q_{d_{0}}(f)$.
By Proposition 1.5, $Q_{d}(f)\not\hookrightarrow Q_{d}$ for $d\geq d_{0}$.
$\hfill\Box$

Lemmas 2.2--2.5 also give several bad strings consisting of 4 blocks.
Every one of those lemmas can be proved by a similar argument in Lemma 2.1.
We omit the detail proof of every lemma, but present the $p$-critical words for $Q_{d_{0}}(f)$.

\trou\noi{\bf Lemma 2.2.} \emph{If $x_{1}=1$ and $y_{n}\geq y_{1}+2$,
then $Q_{d}(f)\not\hookrightarrow Q_{d}$ for
$d\geq 2\sum\limits_{i=2}^{n}(x_{i}+y_{i})+y_{1}$.}

\trou\noi{\bf Proof.} Let $d_{0}=2\sum\limits_{i=2}^{n}(x_{i}+y_{i})+y_{1}$. Then

 $\alpha=10^{y_{1}}1^{x_{2}}0^{y_{2}}\cdots1^{x_{n-1}}0^{y_{n-1}}
 1^{x_{n}}0^{y_{n}-y_{1}-2}00^{y_{1}}11^{x_{2}-1}0^{y_{2}}
 \cdots1^{x_{n-1}}0^{y_{n-1}}1^{x_{n}}0^{y_{n}}$ and

$\beta=10^{y_{1}}1^{x_{2}}0^{y_{2}}\cdots1^{x_{n-1}}0^{y_{n-1}}
 1^{x_{n}}0^{y_{n}-y_{1}-2}10^{y_{1}}01^{x_{2}-1}0^{y_{2}}
\cdots1^{x_{n-1}}0^{y_{n-1}}1^{x_{n}}0^{y_{n}}$ are $2$-critical words for $Q_{d_{0}}(f)$.
By Proposition 1.5, $Q_{d}(f)\not\hookrightarrow Q_{d}$ for $d\geq d_{0}$.
$\hfill\Box$

\trou\noi{\bf Lemma 2.3.}
\emph{If $x_{n}\geq x_{1}+2$ and $y_{1}\geq y_{n}+2$,
then $Q_{d}(f)\not\hookrightarrow Q_{d}$ for
$d\geq 2\sum\limits_{i=2}^{n}x_{i}+2\sum\limits_{i=1}^{n-1}y_{i}+x_{1}+y_{n}-2$.}

\trou\noi{\bf Proof.} Let $d_{0}\geq2\sum\limits_{i=2}^{n}x_{i}+2\sum\limits_{i=1}^{n-1}y_{i}+x_{1}+y_{n}-2$. Then

$\alpha=1^{x_{1}}0^{y_{1}}1^{x_{2}}0^{y_{2}}\cdots1^{x_{n-1}}0^{y_{n-1}}1^{x_{n}-2}
010^{y_{1}-2}1^{x_{2}}0^{y_{2}}\cdots1^{x_{n-1}}0^{y_{n-1}}1^{x_{n}}0^{y_{n}}$ and

$\beta=1^{x_{1}}0^{y_{1}}1^{x_{2}}0^{y_{2}}\cdots1^{x_{n-1}}0^{y_{n-1}}1^{x_{n}-2}
100^{y_{1}-2}1^{x_{2}}0^{y_{2}}\cdots1^{x_{n-1}}0^{y_{n-1}}1^{x_{n}}0^{y_{n}}$ are $2$-critical words for $Q_{d_{0}}(f)$.
By Proposition 1.5, $Q_{d}(f)\not\hookrightarrow Q_{d}$ for $d\geq d_{0}$.
$\hfill\Box$

\trou\noi{\bf Lemma 2.4.} \emph{If $x_{n}\geq x_{1}+1$ and $y_{1}=y_{n}$,
then $Q_{d}(f)\not\hookrightarrow Q_{d}$ for
$d\geq 2\sum\limits_{i=2}^{n}x_{i}+2\sum\limits_{i=2}^{n-1}y_{i}+3y_{1}+x_{1}-1$.}

\trou\noi{\bf Proof.}
Let $d_{0}=2\sum\limits_{i=2}^{n}x_{i}+2\sum\limits_{i=2}^{n-1}y_{i}+3y_{1}+x_{1}-1$. Then

$\alpha=1^{x_{1}}0^{y_{1}}1^{x_{2}}0^{y_{2}}\cdots 1^{x_{n}-1}10^{y_{1}-1}11^{x_{2}-1}0^{y_{2}}\cdots
1^{x_{n-1}}0^{y_{n-1}}1^{x_{n}}0^{y_{1}}$ and

$\beta=1^{x_{1}}0^{y_{1}}1^{x_{2}}$
$0^{y_{2}}\cdots1^{x_{n}-1}00^{y_{1}-1}01^{x_{2}-1}0^{y_{2}}\cdots
1^{x_{n-1}}0^{y_{n-1}}1^{x_{n}}0^{y_{1}}$ are $2$-critical words for $Q_{d_{0}}(f)$.
By Proposition 1.5, $Q_{d}(f)\not\hookrightarrow Q_{d}$ for $d\geq d_{0}$.
$\hfill\Box$

\trou\noi{\bf Lemma 2.5.} \emph{If $y_{1}\geq y_{n}+1$ and $x_{1}=x_{n}=1$,
then $Q_{d}(f)\not\hookrightarrow Q_{d}$ for
$d\geq 2\sum\limits_{i=2}^{n-1}x_{i}+2\sum\limits_{i=1}^{n-1}y_{i}+y_{n}+2$.}

\trou\noi{\bf Proof.}
Let $d_{0}=2\sum\limits_{i=2}^{n-1}x_{i}+2\sum\limits_{i=1}^{n-1}y_{i}+y_{n}+2$. Then

 $\alpha=10^{y_{1}}1^{x_{2}}0^{y_{2}}\cdots1^{x_{n-1}}0^{y_{n-1}-1}000^{y_{1}-1}
1^{x_{2}}0^{y_{2}}\cdots1^{x_{n-1}}0^{y_{n-1}}10^{y_{n}}$ and

$\beta=10^{y_{1}}1^{x_{2}}$
$0^{y_{2}}\cdots1^{x_{n-1}}0^{y_{n-1}-1}110^{y_{1}-1}
1^{x_{2}}0^{y_{2}}\cdots1^{x_{n-1}}0^{y_{n-1}}10^{y_{n}}$ are $2$-critical words for $Q_{d_{0}}(f)$.
By Proposition 1.5, $Q_{d}(f)\not\hookrightarrow Q_{d}$ for $d\geq d_{0}$.
$\hfill\Box$

If $x_{1}\geq2$ and $y_{n}\geq2$, then $f$ is bad by Lemma 2.1.
So $f$ is good only if one of $x_{1}$ and $y_{n}$ is 1.
Without loss of generality let $x_{1}=1$ and $y_{n}\geq1$ by Proposition 1.3.

 \trou\noi{\bf Proof of Theorem 1.1.}
Obviously, the good strings can be found in those which not be mentioned in Lemmas 2.1-2.6.
As the above discussion, $f$ is good only if $x_{1}=1$ and $y_{n}\geq1$.
If $x_{1}=1$ and $y_{n}\geq y_{1}+2$, then $f$ is bad by Lemma 2.2.
So $f$ is good only if $x_{1}=1$ and $y_{n}\leq y_{1}+1$.
If $x_{1}=1$, $x_{n}\geq3$ and $y_{1}\geq y_{n}+2$, then $f$ is bad by Lemma 2.3.
So $f$ is good only if
$x_{1}=1$, $y_{n}= y_{1}+1$ and $x_{n}\geq 1$;
$x_{1}=1$, $y_{n}= y_{1}$ and $x_{n}\geq 1$;
$x_{1}=1$, $y_{n}+1= y_{1}$ and $x_{n}\geq 1$;
or $x_{1}=1$, $y_{1}\geq y_{n}+2$ and $2\geq x_{n}\geq 1$.
If $x_{1}=1$, $x_{n}\geq2$ and $y_{n}= y_{1}$, then $f$ is bad by Lemma 2.4.
So $f$ is good only if $x_{1}=1$, $y_{n}= y_{1}+1$ and $x_{n}\geq 1$;
$x_{1}=1$, $y_{n}= y_{1}-1$ and $x_{n}\geq 1$;
$x_{1}=1$, $y_{n}= y_{1}$ and $x_{n}=1$;
or $x_{1}=1$, $y_{1}\geq  y_{n}+2$ and $2 \geq x_{n}\geq 1$.
If $x_{1}=1$, $x_{n}=1$ and $y_{1}\geq y_{n}+1$, then $f$ is bad by Lemma 2.5.
So if good only if (a) $x_{1}=1$ and $y_{n}=y_{1}+1$;
(b) $x_{1}=x_{n}=1$ and $y_{1}=y_{n}$;
(c) $x_{1}=1$, $x_{n}\geq2$ and $y_{1}=y_{n}+1$;
or (d) $x_{1}=1$, $x_{n}=2$ and $y_{1}\geq y_{n}+2$.
$\hfill\Box$

\section{The index of string consisting of 4 blocks}
To prove Theorem 1.2, we need show that $Q_{d}(f)\hookrightarrow Q_{d}$ for every string $f$ in Table 2 if and only if $d<B(f)$, and Table 2 covers $F$.
For the former, we prove it by Lemmas 3.1 and 3.2.
For the latter, we prove it by Lemma 3.3.

\trou\noi{\bf Lemma 3.1.}
\emph{If $f$ is one of strings $(1)$-$(4)$, then it is good.}

\trou\noi{\bf Proof} Let $n=2$, $x_{1}=r$, $y_{1}=s$, $x_{2}=t$ and $y_{2}=k$.
Then all the good strings from $F$ must among $(1)$-$(4)$ by Theorem 1.1.
For the string $(2)$, it is good by Theorems 4.3 and 4.4 in $\cite{asy}$.
To save space we only show that $(4)$ is good indeed since $(1)$ and $(3)$ can be prowed to be good by similar method.

Let $d\geq1$, $\alpha=a_{1}a_{2}\cdots a_{d}$,
$\beta=b_{1}b_{2}\cdots b_{d} \in V(Q_{d}(f))$,
$a_{i_{j}}\neq b_{i_{j}}$ and $\alpha_{j}=\alpha+e_{j}$,
where $j=1,\ldots,p$.
Without loss of generality, assume $a_{i_{1}}=0$.
Then the following claim holds.

\trou\noi{\bf Claim:} \emph{If $p\geq2$, then there exists $k\in\{1,2,\ldots,p\}$ such that $\alpha_{k}\in V(Q_{d}(f))$}.

\noindent
\emph{Proof.}
On the contrary we suppose that $\alpha_{k}\notin V(Q_{d}(f))$ for all $k\in\{1,2,\ldots,p\}$.
Then we can show that $p\geq 3$ by all the three possible cases of the factor of $\alpha_{1}$ containing $a_{i_{1}}$.

\textbf{Case 1.}
$\alpha_{i_{1},i_{1}+s+k+2}=00^{s}1^{2}0^{k}$.

Since $\beta \in V(Q_{d}(f))$ and $\alpha_{2}\not\in V(Q_{d}(f))$, $i_{2}\in[i_{1}+1,i_{1}+s+k+2]$.
As $s\geq k+2$, there are two following cases for the coordinate of $i_{2}$.

\textbf{Subcase 1.1.} $i_{2}=i_{1}+s+2$.

In this subcase $\alpha_{i_{1},i_{2}+s+k+1}=00^{s}110^{s-1}1^{2}0^{k}$.
If $p=2$, then $\beta_{i_{2}-1,i_{2}+s+k+1}=100^{s-1}1^{2}0^{k}=f$,
a contradiction since $\beta \in V(Q_{d}(f))$.
So $p\geq3$.

\textbf{Subcase 1.2.} $i_{2}=i_{1}+s+t+2$ for some $t$ such that $1\leq t \leq k$.

In this subcase $\alpha_{i_{1},i_{2}+s+k+3}=00^{s}110^{t-1}00^{s}1^{2}0^{k}$.
If $p=2$, then $\beta_{i_{2},i_{2}+s+k+2}=10^{s}1^{2}0^{k}=f$,
a contradiction since $\beta \in V(Q_{d}(f))$.
So $p\geq3$.

\textbf{Case 2.}
$a_{i_{1}}$ is preceded by $\alpha_{i_{1}-s-1,i_{1}+k+1}=10^{s}010^{k}$.

Since $\beta \in V(Q_{d}(f))$, $\alpha_{2}\not\in V(Q_{d}(f))$ and $s\geq k+2$,
there exists some $t\in [1,k]$ such that $i_{2}=i_{1}+t+1$.
Hence $\alpha_{i_{1}-s-1,i_{2}+s+k+2}=10^{s}010^{t-1}00^{s}1^{2}0^{k}$.
If $p=2$, then $\beta_{i_{2},i_{2}+s+k+2}=10^{s}1^{2}0^{k}=f$,
a contradiction since $\beta \in V(Q_{d}(f))$.
So $p\geq 3$.

\textbf{Case 3.}
$a_{i_{1}}$ is preceded by $\alpha_{i_{1}-s-2,i_{1}+k}=10^{s}100^{k}$.

Since $\beta \in V(Q_{d}(f))$, $\alpha_{2}\not\in V(Q_{d}(f))$ and $s\geq k+2$,
there exists some $t\in [1,k]$ such that $i_{2}=i_{1}+t$.
So $\alpha_{i_{1}-s-2,i_{2}+s+k+2}=10^{s}100^{t-1}00^{s}1^{2}0^{k}$.
If $p=2$, then $\beta_{i_{2},i_{2}+s+k+2}=10^{s}1^{2}0^{k}=f$,
a contradiction since $\beta \in V(Q_{d}(f))$.
So $p\geq 3$.

By the above discussions, we know that $p\geq 3$ and for $j=2$ it satisfies that either
(A) $\alpha_{i_{j}-1,i_{j}+s+k+1}=110^{s-1}110^{k}$ or
(B) $\alpha_{i_{j},i_{j}+s+k+2}=00^{s}110^{k}$.
Now we prove that it also satisfies for $j\geq 3$ by induction on $j$.
We assume that it holds for $j$ such that $2\leq j<d$,
now we prove it holds for $j+1$ under the inductive assumption (A) and (B), respectively.

\textbf{(A).} $\alpha_{i_{j}-1,i_{j}+s+k+1}=110^{s-1}110^{k}$.

As $\beta \in V(Q_{d}(f))$, $\alpha_{j+1}\not\in V(Q_{d}(f))$ and $s\geq k+2$,
there are two possible cases of $i_{j+1}$.

(A.1) $a_{i_{j+1}}=1$ and $i_{j+1}=i_{j}+s+1$.

In this subcase, $a_{i_{j+1}}$ is preceded by $\alpha_{1,i_{j}-2}110^{s-1}1$ and followed by $0^{s-1}1^{2}0^{k}\alpha_{i_{j}+s+k+2,d}$.
So $\alpha_{i_{j+1}-1,i_{j+1}+s+k+1}=110^{s-1}110^{k}$.

(A.2) $a_{i_{j+1}}=0$ and $i_{j+1}=i_{j}+s+t'+1$ for some $t'\in [1,k]$.

In this case, $a_{i_{j+1}}$ is preceded by $\alpha_{1,i_{j}-2}110^{s-1}110^{t'-1}$ and followed by $0^{s}1^{2}0^{k}\alpha_{i_{j+1}+s+k+3,d}$.
So $\alpha_{i_{j+1},i_{j+1}+s+k+2}=00^{s}110^{k}$.

\textbf{(B).} $\alpha_{i_{j},i_{j}+s+k+2}=00^{s}110^{k}$.

Since $\beta \in V(Q_{d}(f))$ and $\alpha_{j+1}\not\in V(Q_{d}(f))$,
there are two cases for the position of $i_{j+1}$.

(B.1) $a_{i_{j+1}}=1$ and $i_{j+1}=i_{j}+s+2$.

In this case, $a_{i_{j+1}}$ is preceded by $\alpha_{1,i_{j}-1}00^{s}1$ and followed by $0^{s-1}1^{2}0^{k}\alpha_{i_{j}+s+k+2,d}$.
So $\alpha_{i_{j+1}-1,i_{j+1}+s+k+1}=110^{s-1}110^{k}$.

(B.2) $a_{i_{j+1}}=0$ and $i_{j+1}=i_{j}+s+t'+2$ for some $t'\in [1,k]$.

In this case, $a_{i_{j+1}}$ is preceded by $\alpha_{1,i_{j}-1}00^{s}110^{t'-1}$ and followed by $0^{s}1^{2}0^{k}\alpha_{i_{j+1}+s+k+3,d}$.
So $\alpha_{i_{j+1},i_{j+1}+s+k+2}=00^{s}110^{k}$.

Thus, $\alpha_{i_{j+1},i_{j+1}+s+k+2}=00^{s}110^{k}$
or $\alpha_{i_{j+1}-1,i_{j+1}+s+k+1}=110^{s-1}110^{k}$, that is, (A) or (B) holds for $j+1$.
Specially, $\alpha_{i_{p},i_{p}+s+k+2}=00^{s}110^{k}$
or $\alpha_{i_{p}-1,i_{p}+s+k+1}=110^{s-1}110^{k}$.
Hence $\beta_{i_{p},i_{p}+s+k+2}=10^{s}110^{k}=f$
or $\beta_{i_{p}-1,i_{p}+s+k+1}=100^{s-1}110^{k}=f$.
It is a contradiction since $\beta\in V(Q_{d}(f))$.
Thus the claim holds.

Now we show that $H(\alpha,\beta)=d_{Q_{d}(f)}(\alpha,\beta)=p$ by induction on $p$. Obviously, it is trivial for $p=1$.
Suppose $p\geq2$ and it holds for $p-1$.
There exists some $\alpha_{j} \in V(Q_{d}(f))$ by the above claim,
so $p-1=d_{Q_{d}(f)}(\alpha_{j},\beta)$ by the induction hypothesis,
hence $p=d_{Q_{d}(f)}(\alpha,\beta)$ and so $f$ is good. $\hfill\Box$

\trou\noi{\bf Lemma 3.2.}
\emph{If $f$ is one of strings $(5)$-$(14)$, then $B(f)$ is given just as shown in Table 2.}

\trou\noi{\bf Proof.}
Let $f$ be one of strings $(5)$-$(14)$.
We need show that $Q_{d}(f)\not\hookrightarrow Q_{d}$ for $d\geq B(f)$,
and $Q_{d}(f)\hookrightarrow Q_{d}$ for $ d<B(f)$.
For the former, we show there exist $p$-critical words for $Q_{B(f)}(f)$ in the following,
hence we know it holds by Proposition 1.4.

(5): $\alpha=1^{r}0^{s}1^{t}0 1 1^{r-t-2}0^{s}1^{t}0^{2}$
 and $\beta=1^{r}0^{s}1^{t}10 1^{r-t-2}0^{s}1^{t}0^{2}$;

(6): $\alpha=1^{t+2}0^{s-2}011^{t}0^{s}1^{t}0^{k}$ and $\beta=1^{t+2}0^{s-2}101^{t}0^{s}1^{t}0^{k}$;

(7): $\alpha=1^{r}0^{s-1}01^{t}11^{r-t-2}0^{s}1^{t}0^{k}$ and $\beta=1^{r}0^{s-1}11^{t}01^{r-t-2}0^{s}1^{t}0^{k}$;

(8): $\alpha=1^{r}0^{s-1}01^{t}0^{s}11^{t-1}0^{k}$ and $\beta=1^{r}0^{s-1}11^{t}0^{s}01^{t-1}0^{k}$;

(9): $\alpha=10^{s}100^{s}10^{k}$ and $\beta=10^{s}110^{s}00^{k}$;

(10): $\alpha=1 0^{k}00^{s-k-1}00^{k}10^{k}$ and $\beta=1 0^{k}10^{s-k-1}10^{k}10^{k}$;

(11): $\alpha=1^{r}0^{s-1}01^{r-1}00^{s-1}1^{r}0^{k}$ and $\beta=1^{r}0^{s-1}11^{r-1}10^{s-1}1^{r}0^{k}$;

(12): $\alpha=1^{r}0^{s}1^{t-2}010^{s-2}1^{t}0^{k}$ and $\beta=1^{r}0^{s}1^{t-2}100^{s-2}1^{t}0^{k}$;

(13): $\alpha=1^{r}0^{s}1^{2}0^{s}100^{k}$ and
$\beta=1^{r}0^{s}1^{2}0^{s}010^{k}$; and


(14): $\alpha=1^{r}0^{s}1^{t}0^{k-2}101^{r-2}0^{s}1^{t}0^{k}$ and $\beta=1^{r}0^{s}1^{t}0^{k-2}011^{r-2}0^{s}1^{t}0^{k}$;

(15): $\alpha=1^{r}0^{s}1^{t}010^{s}1^{t}0^{k}$ and $\beta=1^{r}0^{s}1^{t}100^{s}1^{t}0^{k}$.

For the latter, to save space here we only give the proof of $(12)$,
the others can be proved by similar discussions.
Let $ d<B(f)$. Assume $\alpha=a_{1}a_{2}\cdots a_{d}$ and $\beta=b_{1}b_{2}\cdots b_{d}$
be any vertices of  $Q_{d}(f)$, $a_{i_{j}}\neq b_{i_{j}}$ and $\alpha_{j}=\alpha+e_{j}$, where $j=1,\ldots,p$.
Without loss of generality, suppose that $a_{i_{1}}=1$. Then we have the following claim:

\trou\noi{\bf Claim} \emph{Let $p\geq2$. Then there exists $k\in\{1,2,\ldots,p\}$ such that $\alpha_{k}\in V(Q_{d}(f))$}.

\noindent
\emph{Proof.}
On the contrary we suppose that $\alpha_{k}\notin V(Q_{d}(f))$ for all $k\in\{1,2,\ldots,p\}$.
Then we can show that $p=1$ by the following cases.
Hence we get a contradiction since $p\geq2$.

\textbf{Case 1.} $\alpha_{i_{1}-r-s_{1},i_{1}+s_{2}+t+k}=1^{r}0^{s_{1}}10^{s_{2}}1^{t}0^{k}$,
where $s_{1}+s_{2}=s-1$.

Obviously by $\beta \in V(Q_{d}(f))$, $i_{2}\in [i_{1}+1,i_{1}+s_{2}+t+k]$.
Since $\alpha_{2}\not\in V(Q_{d}(f))$, $t\geq r+2$ and $s\geq k+2$,
the position of $i_{2}$ may be the following three subcases.

\textbf{Subcase 1.1.} $i_{2}=i_{1}+1$, $s_{1}=0$ and $s_{2}=s-1$.

In this subcase $\alpha_{i_{1}-r-s-t+2,i_{2}+s+t+k-2}=1^{r}0^{s}1^{t-2}100^{s-2}1^{t}0^{k}$.
But $d=|\alpha|\geq 2s+2t+r+k-2=B(f)$, a contradiction.

\textbf{Subcase 1.2.} $i_{2}=i_{1}+s_{2}+t$.

In this case $\alpha_{i_{1}-r-s_{1},i_{2}+s+k+t-1}= 1^{r}0^{s_{1}}10^{s_{2}}1^{t-1}10^{s-1}1^{t}0^{k}$.
But $d=|\alpha|\geq 2s+2t+r+k-1>B(f)$, a contradiction.

\textbf{Subcase 1.3.} $i_{2}=i_{1}+s_{2}+t+k$.

In this case $\alpha_{i_{1}-r-s_{1},i_{2}+r+s+k+t-1}
=1^{r}0^{s_{1}}10^{s_{2}}1^{t}0^{k-1}01^{r-1}0^{s}1^{t}0^{k} $.
But $d=|\alpha|\geq 2s+2t+2r+2k-1>B(f)$, a contradiction.

\textbf{Case 2.} $\alpha_{i_{1}-r-s-t-k_{1},i_{1}+k_{2}}=1^{r}0^{s}1^{t}0^{k_{1}}10^{k_{2}}$, where $k_{1}+k_{2}=k-1$, $k_{1}\geq 0$ and $k_{2}\geq 1$.

We distinguish two subcases by $r=1$ and $r\geq 2$ to discussion.

\textbf{Subcase 2.1.} $r=1$.

In this case, $i_{2}=i_{1}+k_{3}+1$ for some $k_{3}\in [0,k_{2}-1]$ and
$\alpha_{i_{1}-r-s-t-k_{1},i_{2}+s+t+k}
=1^{r}0^{s}1^{t}0^{k_{1}}10^{k_{3}}00^{s}1^{t}0^{k}$.
But $d=|\alpha|\geq 2s+2t+r+k+2+k_{1}+k_{3}>B(f)$, a contradiction.

\textbf{Subcase 2.2.} $r\geq 2$.

In this subcase, there are two possible positions for $i_{2}$.

One is that $i_{2}=i_{1}+k_{2}$ and
$\alpha_{i_{1}-r-s-t-k_{1},i_{2}+r+s+t+k-1}
=1^{r}0^{s}1^{t}0^{k_{1}}10^{k_{2}-1}01^{r-1}0^{s}1^{t}0^{k}$.
But $|\alpha|\geq 2r+2s+2t+2k-1>B(f)$,
a contradiction.
The other one is that $i_{2}=i_{1}+1$, $k_{2}=1$ and
$\alpha_{i_{1}-r-s-t-k+2,i_{2}+r+s+t+k-2}$
$=1^{r}0^{s}1^{t}0^{k_{1}}10^{k_{2}-1}01^{r-1}0^{s}1^{t}0^{k}$.
But $|\alpha|\geq 2r+2s+2t+2k-2>B(f)$, a contradiction.

Hence by the above discussion, if $\alpha_{k}\in V(Q_{d}(f))$ for all $k\in\{1,2,\ldots,p\}$,
then there exist no $i_{2}$, that is, $p=1$, a contradiction to $p\geq2$.
Thus the claim holds.

Now we show that $H(\alpha,\beta)=d_{Q_{d}(f)}(\alpha,\beta)=p$ by induction on $p$,
where $ d<B(f)$.
Obviously it holds for $p=1$.
Suppose $p\geq 2$ and it holds for $p-1$.
By the above claim there exists $\alpha_{j}\in  V(Q_{d}(f))$,
so $H(\alpha_{j},\beta)=d_{Q_{d}(f)}(\alpha_{j},\beta)=p-1$ by the induction hypothesis.
Hence $H(\alpha,\beta)=d_{Q_{d}(f)}(\alpha,\beta)=p$, in other words,
$Q_{d}(f)\hookrightarrow Q_{d}$ for $ d<B(f)$.

\trou\noi{\bf Lemma 3.3.} \emph{Table 2 covers $F$.}

\noindent
\emph{Proof.}
First we show a relation between $r$, $s$, $t$ and $k$ by the following 36 cases in Table 3,
where all of $r$, $s$, $t$ and $k$ are positive integers.
It is not difficult to find that $(i.j)$ and $(j.i)$ are the same case by Proposition 1.3.
Thus we only need to discuss the cases $(i.j)$ that $i\geq j$, where $i,j=1,2,\ldots,6$.

{\renewcommand\baselinestretch{0.9}\selectfont
\begin{center} \fontsize{10.5pt}{11.5pt}\selectfont
\begin{longtable}{l|l|l|l|l|l|l}
\caption{\label{comparison} A relation between $r$, $s$, $t$ and $k$.}\\
\hline
&$t\geq r+2$&$t=r+1$&$t=r$&$r=t+1$&$r=t+2$&$r\geq t+3$ \\
\hline
$s\geq k+2$&(1.1) &(2.1) &(3.1) &(4.1)& (5.1)& (6.1)\\
\hline
 $s=k+1$ &(1.2) &(2.2) &(3.2) &(4.2)& (5.2)& (6.2) \\
\hline
$s=k$&(1.3) &(2.3) &(3.3) &(4.3)& (5.3)& (6.3) \\
\hline
$k=s+1$ &(1.4) &(2.4) &(3.4) &(4.4)& (5.4)& (6.4)\\
\hline
 $k=s+2$  &(1.5) &(2.5) &(3.5) &(4.5)& (5.5)& (6.5) \\
\hline
 $k\geq s+3$ &(1.6) &(2.6) &(3.6) &(4.6)& (5.6)& (6.6)\\
\hline
\end{longtable}
\end{center}
\par}

\noindent
 {\bf Case (1.1)}: $t\geq r+2$ and $s\geq k+2$.

This case is covered by $(12)$.

\noindent
 {\bf Cases (2.1) and (2.2)}: $t= r+1$ and $s\geq k+1$.

For the subcase $r\geq 2$ and $k\geq 2$, it is covered by $(15)$.
For the subcase $r\geq 1$ and $k=1$, its reverse is covered  by $(3)$.
For the subcase $r=1$, $k\geq2$ and $s=k+1$, it is covered by $(3)$.
For the subcase $r=1$, $k\geq2$ and $s\geq k+2$, it is covered by $(4)$.

\noindent
 {\bf  Cases (3.1) and (3.2)}: $t= r$ and $s\geq k+1$.

For the subcases $r=1$ and $r\geq2$ it is covered by $(10)$ and $(11)$, respectively.

\noindent
 {\bf   Case (3.3)}: $t= r$ and $s=k$.

If $r\geq 2$ and $k\geq 2$, then it is covered by $(15)$,
and $r=1$ and $k\geq 1$  is covered by $(2)$.

\noindent
 {\bf Cases (4.1), (4.2) and (4.3)} : $r= t+1$ and $s\geq k$.

If $k=1$, then its reverse is covered by $(1)$.
If $k\geq 2$, then it is covered by $(15)$.

\noindent
 {\bf  Cases (4.4)}: $r= t+1$ and $k= s+1$.

This case is covered by $(8)$.

\noindent
 {\bf  Cases (5.1), (5.2) and (5.3)}: $r= t+2$ and $s\geq k$.

If  $s=1$, then it is covered by $(7)$.
If $s\geq 2$, then it is covered by $(6)$.

\noindent
 {\bf Cases (5.4)}: $r= t+2$ and $k=s+1$.

If  $s=1$, then its reverse is covered by $(8)$.
If $s\geq 2$, then it is covered by $(7)$.

\noindent
 {\bf Cases (5.5)}: $r= t+2$ and $k=s+2$.

If $s=1$ and $t=1$, or $s\geq 2$ and $t\geq1$, then it is covered by $(8)$ or $(15)$, respectively.

\noindent
 {\bf Cases (6.1) and (6.2)}: $r\geq t+3$ and $s\geq k+1$.

If $s=2$, then it is covered by $(13)$.
If $s\geq 3$ and $k=1$, then it is covered by $(7)$.
If $s\geq 3$ and $k=2$, then it is covered by $(5)$.
If $s\geq 3$ and $k\geq 3$, then it is covered by $(14)$.

\noindent
 {\bf Cases (6.3) }: $r\geq t+3$ and $k=s$.

The reverse of subcases $k=1$, $k\geq 2$ and $k\geq 3$ is covered by $(9)$, $(13)$ and $(14)$, respectively.

\noindent
 {\bf  Cases (6.4)}: $r\geq t+3$ and $k=s+1$.

If  $s=1$, then its reverse is covered by $(8)$.
If $s\geq 2$, then its reverse  is covered by $(14)$.

\noindent
 {\bf  Cases (6.5)}: $r\geq t+3$ and $k=s+2$.

If $t=1$ and $s=1$, then its reverse is covered by $(8)$.
If $t=1$ and $s\geq 2$, or $t\geq1$ and $s\geq 2$, then its reverse is covered by $(14)$.

\noindent
 {\bf  Cases (6.6)}: $r\geq t+3$ and $k\geq s+3$.

If $s=1$ and $t=1$, then it is covered by $(8)$;
If $s=1$ and $t\geq 2$, or $s\geq 2$ and $t\geq2$, then it is covered by $(14)$.
This completes the proof.
$\hfill\Box$


\end{document}